\documentclass[12pt]{article}
\usepackage{latexsym}
\usepackage{amssymb}
\usepackage{epsfig}

\newcommand{\ds}{\displaystyle}
\newcommand{\beq}{\begin{equation}}
\newcommand{\eeq}{\end{equation}}
\newcommand{\beqr}{\begin{eqnarray}}
\newcommand{\eeqr}{\end{eqnarray}}

\newcommand{\ve}{\varepsilon}

\begin{document}

\title{Design Flaws in the Implementation of the Ziggurat and Monty Python methods
(and some remarks on Matlab randn) }
\author{
Boaz Nadler\footnote{Department of Computer Science and Applied
Mathematics, Weizmann Institute of Science, PO Box 26, Rehovot, Israel.\newline
email address: boaz.nadler@weizmann.ac.il, http://www.wisdom.weizmann.ac.il/$\sim$nadler}}
\maketitle

\begin{abstract}

{\em Ziggurat} and {\em Monty Python} are two fast and elegant
methods proposed by Marsaglia and Tsang to transform uniform
random variables to random variables with normal, exponential and
other common probability distributions. While the proposed methods
are theoretically correct, we show that there are various design
flaws in the uniform pseudo random number generators (PRNG's) of
their published implementations for both the normal and Gamma
distributions \cite{Ziggurat,{Gamma},Monty}. These flaws lead to
non-uniformity of the resulting pseudo-random numbers and
consequently to noticeable deviations of their outputs from the
required distributions. In addition, we show that the underlying
uniform PRNG of the published implementation of Matlab's
\texttt{randn}, which is also based on the Ziggurat method, is not
uniformly distributed with correlations between consecutive pairs.
Also, we show that the simple linear initialization of the
registers in matlab's \texttt{randn} may lead to non-trivial
correlations between output sequences initialized with different
(related or even random unrelated) seeds. These, in turn, may lead
to erroneous results for stochastic simulations.

\end{abstract}

\section{Introduction}

Pseudo random number generators (PRNG) are a key component of
stochastic simulations. Most PRNG's produce sequences of
(seemingly) uniformly distributed real numbers in the
interval $[0,1)$, typically by quantizing a sequence of integers
in the set $\Omega_k = \{0,1,2^k-1\}$ via division by $2^k$.
Random variables with normal (Gaussian), Gamma or other
distributions, are then typically constructed via transformations
of uniform random variables \cite{Knuth,{Devroye}}.

In \cite{Marsaglia84}, Marsaglia and Tsang proposed the Ziggurat
method to generate random variables with a decreasing or symmetric
unimodal density from uniform random variables. While in the
original paper \cite{Marsaglia84} the authors assumed the
existence  of a suitable PRNG that outputs uniform random numbers
in the range $[0,1)$, in \cite{Ziggurat} and \cite{Gamma}, the
same authors refined their original method and suggested a
specific underlying fast uniform PRNG along with computer code to
produce normal or Gamma distributed numbers. In a different work
\cite{Monty}, the same authors proposed the Monty Python method to
generate random variables with normal, exponential or other common
distributions, and proposed a
specific implementation based on a multiply-with-carry (MWC)
uniform PRNG. In this paper we show that there are various flaws
in the design of these underlying PRNG's, which lead to
significant deviations of their outputs from uniformity, and consequently poor distributions of the
resulting normal or Gamma distributions.

We note that statistical problems with the implementation of \cite{Ziggurat}
were recently noted by Leong \& al \cite{Leong}. These authors found that the resulting sequences of normally
distributed numbers fail the simple $\chi^2$ test, and attributed this finding to the
(relatively) short period of the underlying PRNG ($2^{32}-1$). Our analysis elucidates the
design flaws leading to these statistical problems, which are mainly due to
non-uniformity and correlations of the outputs of the PRNG.
While the short period of the specific suggestion of \cite{Ziggurat,Gamma},
only magnifies these problems, we show that other PRNG's with much longer periods but same
output function lead to the same statistical problems.

Another normal random number generator based on the Ziggurat
method is Matlab's built-in function {\tt randn}. We analyze
the underlying uniform PRNG of this
function, based on the matlab code published in
\cite{matlab_randn}. We show that while individual outputs of this
PRNG are uniformly distributed, pairs of
consecutive outputs are correlated. Since the output of the Ziggurat method is highly
non-linear in its input, it is difficult to detect these
correlations in the resulting normal random numbers. However,
we show that initializations of the function \texttt{randn} with
different, either related or even random unrelated seeds,
as done in parallel implementations and
other stochastic simulations, may lead
to non-trivial correlations between the resulting output sequences
of random numbers. We give a simple example where a sequence of
such initializations yields incorrect results for a stochastic
simulation.

The paper is organized as follows. In section 2 we present a
probabilistic setting in which we can properly define statistical
properties of sequences of random numbers from PRNG's. The design
flaws and statistical weaknesses of the PRNG's published in
\cite{Ziggurat,Gamma} are analyzed in section 3. In Section 4 we
analyze the weaknesses in the multiply with carry
generator and their consequences on the Monty Python method of
\cite{Monty}. The analysis of Matlab's {\tt randn} is presented in
Section 5. We conclude in Section 6 with a summary and
discussion.

\section{Randomness and statistical requirements from a PRNG  }

The main goal of a uniform PRNG on a set of possible outputs $\Omega$
is to produce long sequences of numbers which imitate, in a statistical sense,
realizations of a sequence of corresponding independently identically distributed (i.i.d.)
uniform random variables
on the same set.
While sequences of random variables are a well defined concept within the framework of probability theory,
the notion of "randomness" in a sequence of numbers
produced by a deterministic algorithm is problematic and requires proper definition.

Following \cite{Lecuyer}, we consider a uniform random number
generator as a finite state machine $(F,g,s,P)$ with an internal
state set $\cal S$ and output set $\Omega$. The current state of
the PRNG is denoted by $s\in{\cal S}$, $F:{\cal S}\to{\cal S}$ is
the transition function, $g:{\cal S}\to\Omega$ is the output
function, and $P:{\cal T}\to{\cal S}$ is the initialization
function, where ${\cal T}$ is a set of possible seeds. The
internal state is updated according to $s_{i+1} = F(s_i)$ while
the output at step $i$ is $g(s_i)$. The state machine is
initialized by a seed, such that $s_0 = P(seed)$. We denote by
$u(s_0,j)$ the $j$-th output of a generator with initial state
$s_0$. Note that $\cal S$ and $\Omega$ need not (and in general
should not) be of the same size. For example, the internal state
could be of size $|{\cal S}| = 2^{128}$ (e.g. 128 bits), while the
output set $\Omega$ could be of size $2^{32}$.

Since the PRNG is deterministic, the output sequence is uniquely determined
by the initial seed and the functions $F,g,P$. Moreover, since
$F$ is a finite state machine, it is obviously
periodic with some period $l\leq |{\cal S}|$.
We introduce a probability
space in this setting by considering both the time of observation
and the initial state $s_0$ as {\em random variables}, with
the initial state  uniformly distributed over the set $\cal S$
and the time of observation uniformly distributed over the integer set $\{1,2\ldots l\}$.
We denote by $\{U_1,U_2,\ldots\}$
the resulting sequence of random variables.

This sequence should, by definition, have the same probability distribution as that of a sequence of
i.i.d. uniform random variables $\{X_1,X_2,\ldots\}$ over the set $\Omega$.
Obviously, the sequence $\{U_i\}$ does not have the same distribution as that of $\{X_i\}$, since
$U_{i+l} = U_i$, for example. Therefore, one of the basic requirements from a PRNG is to have a
very long period $l\gg 1$, significantly longer than the total number of outputs used by the simulation.
In addition, inside the long period we require that at least the {\em low order
statistics} of $\{U_i\}$ and $\{X_i\}$ coincide.
Specifically, we consider a PRNG $(F,g,s,P)$
as statistically sound if it satisfies (at least) the following requirements (see also \cite{Brent}):
\begin{enumerate}

\item {\em Uniformity of first order statistics on $\Omega$:} We
require that
\begin{equation}
\Pr\{U_1 = u\in\Omega\} :\equiv
\frac{1}{|{\cal S}| l }\sum_{s_0\in{\cal S}} \sum_{t=1}^l
\delta(u(s_0,t),u) = \frac{1}{|\Omega|}
    \label{pr_u1}
\end{equation}
where $\delta(i,j)$ is the Kroneker delta function, equal to one if $i=j$ and zero otherwise.

\item {\em Uniformity of second order statistics on $\Omega$:} We
require that
\begin{equation}
  \Pr\{(U_1,U_2)=(u_1,u_2)\} :\equiv \frac{1}{|{\cal S}| l}
  \sum_{s_0\in{\cal S}} \sum_{t=1}^l \delta(u(s_0,t),u_1)
  \delta(u(s_0,t+1),u_2) = \frac1{|\Omega|^2}
    \label{pr_u1_u2}
\end{equation}
Note that combining requirements (\ref{pr_u1}) and (\ref{pr_u1_u2}) implies that
the {\em conditional} distribution of $U_2$ given $U_1$ is also uniform on the set $\Omega$.
In other words, observation of a single output does not affect the distribution of the next output.

\item {\em Insensitivity to initialization with related seeds:}
Recall that the initial state is set according to $s_0=P(seed)$.
Let $U,U'$ be the random variables that correspond to
initializations with seeds that differ by $\Delta$. We require
that for any $\Delta\in{\cal T}\setminus\{0\}$,
\begin{equation}
  \Pr\{U=u,U'=u'\} = \sum_{seed\in{\cal T}} \delta(u(P(seed),j),u)\delta(u(P(seed+\Delta),j),u') = \frac{1}{|\Omega|^2}
  \label{pr_u_utag}
\end{equation}
\end{enumerate}

For many PRNG's, the set of requirements (\ref{pr_u1})-(\ref{pr_u_utag}) does not hold exactly.
We still consider a PRNG as statistically acceptable if the discrepancy between these distributions and
the corresponding uniform ones is extremely small, say below a threshold $\ve$, such that detection
of this discrepancy would require more than $2^{100}$ outputs, for example.

Obviously, PRNG's need to satisfy many more requirements to be considered acceptable from
a statistical point of view, and their output sequences are typically required to pass
various empirical statistical tests (see, for example \cite{Knuth,Brent,DieHard,Hellekalek} and references therein). However, as we shall see below,
each one of the requirements (\ref{pr_u1})-(\ref{pr_u_utag}) is
{\em essential} in the context of both the Ziggurat and the Monty Python methods, and possibly
so in the more general context of stochastic simulations. While requirements (\ref{pr_u1}) and
(\ref{pr_u1_u2}) seem obvious, requirement (\ref{pr_u_utag}) can be quite important when different
runs are made with different seeds, as in parallel implementations of stochastic simulations.

\section{Design Flaws in the uniform RNG of \cite{Ziggurat,{Gamma}}}

For the paper to be reasonably self contained, we briefly describe the basic Ziggurat method.
To generate a non-negative random variable
with a monotonically decreasing density $f(x)$ from a uniform r.v. $U[0,1]$,
we choose $k$ points $0=x_0<x_1<x_2<\ldots<x_{k-1}$, such that
\[
x_i (f(x_{i-1}) - f(x_i)) = x_{k-1} f(x_{k-1}) + \Pr\{x > x_{k-1}\}\quad 1\leq i \leq k-1
\]
Given the set $\{x_i\}_{i=0}^{k-1}$, the Ziggurat method works as follows:
\begin{enumerate}
\item Choose an index $0\leq i<k$ at random with uniform probability $1/k$.
\item Draw a random number $u$ from the uniform distribution $U[0,1]$, and let $x= u x_i$.
If $i\geq 1$ and $x<x_{i-1}$ return $x$.
\item If $i\geq 1$ draw another uniform random number $y$. If $y(f(x_{i-1})-f(x_i))<f(x)-f(x_i)$ return $x$.
\item If $i=0$, generate an $x$ from the tail $x>x_{127}$ and return $x$.
\item Otherwise, return to step 1.
\end{enumerate}
The generation of values from the tail of the distribution is described explicitly below.
In most applications $k$ is chosen to be a power of 2, (typically $k=64,128,256$),
so that choosing a random index with probability $1/k$ is easily done by considering
the first few bits of a random 32 bit uniformly distributed integer, for example.

The key point and the beauty of the Ziggurat method is that if
the two numbers $x$ and $y$ in steps 2 and 3 are indeed
(statistically) random, independent and uniformly distributed over $[0,1)$,
then the output will be a truly normal distributed random variable.
The original publication \cite{Marsaglia84} described the method with an unspecified underlying
uniform PRNG. However, in \cite{Ziggurat,Gamma} the following code for steps 1-3 was suggested
by Marsaglia and Tsang, using $k=128$,
\begin{verbatim}
unsigned long jsr,jz;
long hz,iz;
#define SHR3 (jz=jsr, jsr^=(jsr<<13), jsr^=(jsr>>17),jsr^=(jsr<<5),jz+jsr)
#define RNOR (hz=SHR3, iz=hz&127, (fabs(hz)<kn[iz])? hz*wn[iz] : nfix())
#define UNI (.5 + (signed) SHR3*.2328306e-9)

float nfix()
{   float x,y;
    for(;;){
        x = hz * wn[iz];
        if(iz==0){ // generate an output from the tail
            do{ x=-log(UNI)/x[k-1]; y=-log(UNI);}
            while(y+y<x*x);
            return (hz>0)? r+x : -r-x;
        }
        if(fn[iz]+UNI*(fn[iz-1]-fn[iz]) < exp(-.5*x*x) ) return x;

        hz=SHR3; iz=hz&127; if(fabs(hz)<kn[iz]) return hz*wn[iz];
    }
}
\end{verbatim}
First, a few explanatory words on the code above: The inline code
\texttt{RNOR} produces a normal random number. It first calls
\texttt{SHR3}, which both updates the 32-bit register \texttt{jsr}
and outputs a 32-bit integer, which should be uniformly
distributed in the set $0\ldots 2^{32}-1$. To produce both the
positive and negative parts of the normal distribution, the output
of {\tt SHR3} is assigned to the (signed) variable \texttt{hz} of
type \texttt{long}. The two tables \texttt{kn} and \texttt{wn} are
initialized to store the following values: $\texttt{kn[i]} =
2^{-31}x_{i-1}/x_i$ and $\texttt{wn[i]} = x_i/2^{31}$ for $i\geq
1$ and special values for $i=0$. The procedure \texttt{nfix()}
takes care of steps 3-5, whenever step 2 fails.

The register \texttt{jsr} is updated via a linear transformation made up of three
shifts, hence the name \texttt{SHR3} (shift register 3). For future use we denote this
linear transformation
by \texttt{T}, so that \texttt{jsr}$^{(t+1)}$ = \texttt{T(jsr}$^{(t)}$\texttt{)}.
According to \cite{Marsaglia_T} this transformation has maximal period, e.g. $2^{32}-1$.
Note, however, that the
output of \texttt{SHR3} is \texttt{jsr+T(jsr)} ({\tt mod $2^{32}$}), and not the value of \texttt{jsr} itself. We will
come back to this point later on.

In our analysis we will need estimates on the number of outputs required to
distinguish between a discrete distribution over $k$ values with probabilities $(p_1,\ldots,p_k)$ and one with
probabilities $(q_1,\ldots,q_k)$ with $q_i = p_i(1+\ve_i)$. As shown in the appendix, for a
distinguisher based on the $\chi^2$ statistic, the number of required outputs is of the order of
\begin{equation}
N  = \frac{\sqrt{2k}}{\ds\sum_{j=1}^k p_j \ve_j^2} = O(1/\ve^2)
    \label{N_estimate}
\end{equation}

{\tt Design Flaw \# 1:} The first problem we observe, as also recently noted by Doornik \cite{Doornik},
is that the {\em same} 7 least significant bits of \texttt{hz} are used both for
choosing the random index $0\leq i\leq 127$ in step 1, and for the uniform random number $u$
in step 2 of the algorithm. This obviously introduces some statistical dependencies into the
algorithm in two different locations: The first is that the random numbers from the $i$-th index all end with
the same 7 least significant bits, and the other is in the computation of the rejection probabilities
of step 2. Let us roughly estimate this second deviation and its shortcomings:
For a 32 bit uniform random number and a table
with $k=128$ (e.g. 7 bits), the fact that the last 7 bits are fixed induces errors in the
rejection question at step
2 (whether $u x_i < x_{i-1}$) of the order of $\varepsilon = 1/2^{32-7} = 2^{-25}$ (instead of the quantization
error of $2^{-32}$). Therefore,
according to (\ref{N_estimate}), to detect such a deviation one would need $O(1/\ve^2) = 2^{50}$ outputs. While this design flaw
certainly produces a deviation from the normal distribution, it is quite negligible as compared
to the next design flaw that we now describe.

{\tt Design Flaw \# 2:} The output of \verb"SHR3", of the form
\verb"x+Tx" is highly non-uniform and fails to satisfy the basic
requirement (\ref{pr_u1}). Due to the specific structure of
\verb"T", the output $\texttt{x+Tx}$ is not one-to-one, but rather
a contractive mapping with 1543756180 outputs (about $2^{30.5}$)
not possible at all. Thus the output range of \texttt{SHR3} is
restricted to only about $64\%$ of the possible $2^{32}$ outcomes,
with some values 10 times more probable than expected in a uniform
distribution. Table 1 shows the distribution of the number of
sources of a possible value $y$, e.g. the number of values $x$
such that $x+Tx = y$.

This non-uniform distribution of \texttt{SHR3} yields non-negligible deviations from normality for the resulting normal random numbers.
Due to the structure of {\tt T} one can prove that the lowest seven bits of {\tt x+Tx}
are uniformly distributed.
Therefore, the probability of choosing a specific index $i$ in step 1 of the algorithm is still
$1/128$ as should be. However,
the non-uniformity of the output yields non-negligible deviations in the resulting
variables $u x_i$ and in the expected
rejection probabilities at step 2 of the algorithm.
We now describe the effects of these deviations and estimate the number of outputs
needed to detect them in a simple $\chi^2$ test on the resulting normal numbers.

Let $Z$ denote a standard Gaussian variable with zero mean and unit
variance. For $j=1,\ldots,k-1$ we define $p_j = \Pr\{x_{j-1}<|Z|<x_{j}\}$ and
$p_k = \Pr\{|Z|>x_{k-1}\}$. Let $q_j$ denote the
corresponding probabilities in the Ziggurat algorithm, whose underlying PRNG is {\tt SHR3}.
Then, by definition
\begin{equation}
  q_j = \Pr\{x_{j-1}<|x|<x_{j}\} = \sum_{i=0}^{k-1} \Pr\{\mbox{index chosen is }i\}
  \Pr\{x_{j-1}<|x|<x_{j} | \mbox{ index }i\}
\end{equation}
As discussed above, the 7 least significant bits of {\tt SHR3} are uniformly
distributed, and therefore $\Pr\{\mbox{index chosen is }i\} = 1/k$. However,
the probabilities $\Pr\{x_{j-1}<|x|<x_{j} | \mbox{index }i\}$, deviate from the theoretically
expected ones. To estimate $q_j$ we performed the following calculation:
We passed over all $2^{32}-1$ possible initial values for the register {\tt jsr}, computed the
first output of {\tt RNOR}, and created a histogram of hits into the 128 bins $[x_{i-1},x_i]$ and $[x_{127},\infty)$.
In table \ref{t:bins} we present the eight bins with the largest
deviations (measured as $p_i \ve_i^2$, where $q_i = p_i(1+\ve_i)$).
Applying formula (\ref{N_estimate}), we estimate that after an order of $2^{30}$ outputs,
these deviations from the normal distribution can be detected with a
$\chi^2$ test on these 128 bins.

\begin{figure}[t]
\mbox{
\begin{minipage}[t] {\textwidth}
\begin{center}
\begin{tabular}{c}
\psfig{figure=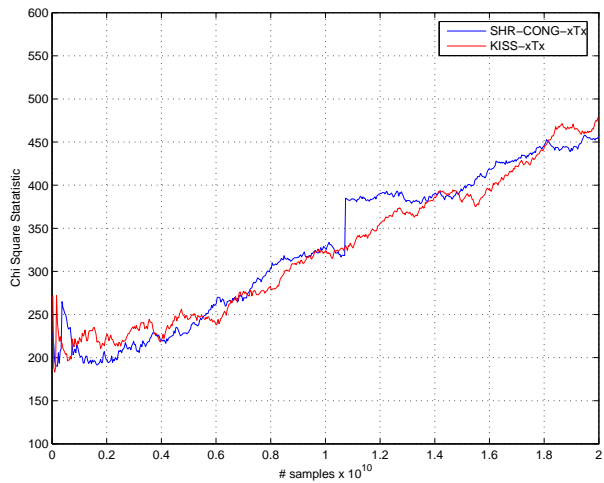,width=8.0cm}\\
\end{tabular}
\end{center}
\end{minipage}
}
\\
\caption{$\chi^2$ values vs. number of samples for generators with
output {\tt x+Tx}. } \label{f:u1}
\end{figure}

This result is not due to the relatively short period of the register {\tt jsr}.
In figure \ref{f:u1}
we present numerical results of the $\chi^2$ test done on 200 bins, evenly spaced in the interval $[-7.0,7.0]$
as done by Leong et al \cite{Leong}, for two other underlying PRNG's with much longer periods,
but whose output is computed via {\tt x+Tx}. The two generators are either a combination
of {\tt SHR} with {\tt CONG}, a multiplicative congruential generator, which is the underlying generator of matlab's {\tt randn}, and the
{\tt KISS} generator \cite{Leong}, which combines also a multiply with carry register.
We stress that in both
cases, the output is {\tt x+Tx} instead of the original {\tt x}, and as expected the $\chi^2$ statistics starts
to significantly deviate from its expected mean after roughly $2^{32}$ outputs.

\begin{table}
\begin{center}
\begin{tabular}{|c|c|}
\hline
\# sources & \# of outputs
\\ \hline
0 & 1543756180 \\ \hline
1 & 1616832933 \\ \hline
2 & 808153149 \\ \hline
3 & 256471123 \\ \hline
4 & 58117590 \\ \hline
5 & 10068341\\ \hline
6 & 1391608\\ \hline
7 & 159565\\ \hline
8 & 15358\\ \hline
9 & 1334\\ \hline
10 & 109\\ \hline
11 & 5\\ \hline
12 & 1\\ \hline
13 & 0\\ \hline
\end{tabular}
\end{center}
\caption{Distribution of the number of sources of the transformation $x+Tx$.}
\end{table}

\begin{table}
\begin{center}
\begin{tabular}{|c|c|c|c|c|c|}
\hline
Interval i & $[x_{i-1},x_i]$ & $p_i$ & $q_i $ & $\ve_i = (q_i-p_i)/p_i$ & $p_i \ve_i^2$
\\ \hline
103    & [2.1443, 2.1659] &    0.0016945  &  0.0016955 &   \ 0.00060306 & 6.1625e-010 \\
82    & [1.7748, 1.7900] &    0.0024789  &  0.0024778   & -0.00041684    & 4.3072e-010 \\
109    & [2.2843, 2.3104] &    0.0014902 &   0.0014894  &  -0.00051947 & 4.0212e-010 \\
92    & [1.9353, 1.9525] &    0.0020838 &   0.0020829  &  -0.00043891 & 4.0143e-010 \\
108   &  [2.2591, 2.2843] &    0.0015242 &   0.001525  & \ 0.00047076 & 3.3779e-010 \\
16    & [0.7981, 0.8189] &    0.0119470 &  0.011945  &  -0.00015320 & 2.8041e-010 \\
104   & [2.1659, 2.1882] &    0.0016603  &  0.001661  & \ 0.00040957 & 2.7852e-010\\
112   &  [2.3659, 2.3954] &    0.001387 &   0.0013864  &  -0.00043609 & 2.6378e-010 \\
\hline
\end{tabular}
\end{center}
\caption{The eight intervals with the largest relative discrepancies
from the normal distribution, measured as $p_i\ve_i^2$.
}\label{t:bins}
\end{table}

\subsection{A quick fix ? Wrong Tail Probabilities !}

Since \verb"x+Tx" is non-uniform, a  possible and natural "quick fix"
is to replace the output by the state of the 32-bit
register \texttt{jsr}, via the following inline code \verb"SHR0",
\begin{verbatim}
#define SHR0 (jsr^=jsr<<13, jsr^= jsr>>17, jsr^=jsr<<5, jsr)
\end{verbatim}
Since \texttt{jsr} is a maximal length 32-bit shift
register, when averaged over its period, its state is a uniformly
distributed 32-bit integer number in the range $1\ldots(2^{32}-1)$, and thus approximately satisfies
requirement (\ref{pr_u1}). One may thus be tempted to conclude that replacing
{\tt SHR3} with {\tt SHR0} fixes all statistical problems, and the resulting sequence of normal
random numbers from the Ziggurat algorithm
should easily pass the $\chi^2$ statistical test of \cite{Leong}.

However, the Ziggurat method with this underlying
PRNG also fails the $\chi^2$ test. The reason is that even though the outputs now satisfy requirement (\ref{pr_u1}),
they fail to satisfy requirement (\ref{pr_u1_u2}), and this leads to non-negligible
deviations in the tail probabilities (when $|x|>x_{127}$).
As seen from the code {\tt RNOR} and the function {\tt nfix()}, a number from the tail is
produced only when the seven least significant bits of {\tt jsr} are
all zero and in addition, the resulting number satisfies the condition {\tt fabs(hz)< kn[0]}.
By enumeration over all $2^{25}-1$ possible values of  \verb"jsr" with 7 least significant bits all zero,
only 2,444,151 values pass this test. For each of these
numbers we calculated the resulting normal number and produced a histogram according
to the following eight intervals defined by $X^t = \{x_{127},3.75,4.0,4.25,4.5,4.75,5.0,5.5\}$,
where $x_{127} = 3.44262$.
In table \ref{t:tail_prob} we show the resulting conditional probabilities $\Pr\{X^t_j<|x|<X^t_{j+1}\,\Big{|}\,|x|>x_{127}\}$
as computed from the Ziggurat method with {\tt SHR0}, vs. the correct probabilities from
the normal distribution. Applying (\ref{N_estimate}) with $k=8$ we obtain that these deviations
can be detected after observation of $O(2^{21})$ samples with $|x|>x_{127}$ or roughly
$N=O(2^{31.5})$ outputs overall.

\begin{table}
\begin{center}
\begin{tabular}{|l|c|c|}
\hline
Interval &$\Pr\{x^t_j<|x|<x^t_{j+1}\}$ & $\Pr\{x^t_j<|Z|<x^t_{j+1}\}$ \\ \hline
  1 &  6.9298e-1&  6.9305e-1\\ \hline
  2 & 1.9705e-1&  1.9700e-1\\ \hline
  3 & 7.2872e-2&  7.2843e-2\\ \hline
  4&  2.5334e-2&  2.5311e-2\\ \hline
  5 & 8.2683e-3&  8.2644e-3\\ \hline
  6 & 2.5105e-3&  2.5357e-3\\ \hline
  7 & 9.3857e-4&  9.2921e-4\\ \hline
  8 & 5.4825e-5&  6.5924e-5\\ \hline
\end{tabular}
\end{center}
\caption{Tail probabilities from the Ziggurat method vs.
the correct probabilities from the normal distribution.
}
\label{t:tail_prob}
\end{table}

\section{Design Flaws in the MWC generator of \cite{Monty}}

In contrast to the Ziggurat method, which requires pre-computation and storage of tables of size $k$,
the Monty Python method \cite{Monty} can produce random variables with normal and other common distributions without
the aid of auxiliary tables.
In \cite{Monty}, the authors presented the method and suggested the following
multiply-with-carry (MWC) generator as a source of uniformly distributed numbers
in the set $\Omega_{32}$,
\begin{verbatim}
unsigned long jsr_z, jsr_w;

#define ZNEW (jsr_z=36969*(jsr_z&65535)+(jsr_z>>16))
#define WNEW (jsr_w=18000*(jsr_w&65535)+(jsr_w>>16))
#define MWC  ((ZNEW<<16) + (WNEW&65535))
\end{verbatim}
In this specific suggestion, each call to {\tt MWC} outputs a 32-bit integer, based on two independent 32-bit registers \verb"jsr_z"
and \verb"jsr_w" . Some of the properties of the transition function of these registers, which is of the form
\begin{equation}
    (c,w) = ((c+a w) \mbox{div }b,(c + a w)\mbox{mod } b)
\end{equation}
where $c$ is the carry and $w$ is the residual upon division by $b$ have been analyzed by
Marsaglia \cite{DieHard}, Couture and
L'Ecuyer \cite{Lecuyer_MWC}. For generalizations to recursions with a similar form, see also Goresky and  Klapper\cite{Klapper}. The main properties
depend on the number $m=ab-1$. As proven in
\cite{DieHard}, when $m$ is a safe prime (both $m$ and $(m-1)/2$ are primes), the period of any starting
value $(c,w)$, other than the trivial fixed points $(0,0)$ and $(a-1,b-1)$ is $(m-1)/2$. Since there are $m+1=ab$
possible values for $(c,w)$, it follows that there are two disjoint orbits with this period.
For the specific values of {\tt MWC} above, the period of each non-trivial orbit of \verb"jsr_w" is
589823999, while that of \verb"jsr_z" is 1211400191. Since both periods are prime numbers
the combined period of {\tt MWC} is their product, or roughly  $2^{59.3}$.

Despite the long period, the main design flaw in {\tt MWC} is the {\em non-uniformity} of its output.
Consider for example the 16 most significant bits of {\tt MWC}, e.g., the lowest 16 bits of the register \verb"jsr_z".
When we consider all possible initializations $(c,w)$, these 16 bits are obviously randomly distributed. However,
{\em within each disjoint orbit}, these 16 bits are non-uniformly distributed.
To illustrate this, we computed the exact probabilities for the seven most significant
bits in a given orbit. In a uniform distribution, the probability to obtain each one of the 128 possible outcomes should
be $p_i = 1/128=0.0078125$. However, as shown in table \ref{t:MSBZ}, some outcomes have non-negligible deviations from this value.
Similarly, the 16 output bits of \verb"jsr_w"
are also not uniformly distributed within each orbit.
Since the union of the two disjoint orbits covers almost all possible $ab$ values for $(c,w)$,
if deviations in one orbit
are of the form $q_i=p_i(1+\ve_i)$, then in
the other orbit the corresponding deviations are approximately $q_i' = p_i(1-\ve_i)$.
Therefore, a union of output sequences belonging to disjoint orbits is approximately uniformly distributed.
This might explain why this generator was reported to pass the DIEHARD tests of Marsaglia.
However, within each orbit, this deviation from uniformity is easily detected with simple statistical tests.
Using formula (\ref{N_estimate}), we obtain that after $2^{28}$
outputs, it can be detected by a $\chi^2$ test with 128 bins. We remark that
on modern computers, $2^{28}$ outputs are typically generated in less than 10 seconds.

To conclude, while the period of {\tt MWC} is larger than $2^{59}$,
its 32-bit output fails the basic requirement (\ref{pr_u1}).
These non-negligible deviations from the uniform distribution also lead to
non-negligible deviations from the normal distribution when applying the Monti Python
method. These can be easily detected by a $\chi^2$ test on 16 bins after $2^{30}$ outputs,
and with even fewer outputs if the number of bins is increased.

\begin{table}
\begin{center}
\begin{tabular}{|c|c|c|}
\hline
7 MSB's $i$ &$\Pr\{{\tt (z\_new>>9)}=i\}$ & $\ve_i$ \\ \hline
22 &      0.007818141  &   \ 7.2200e-3 \\ \hline
107 &     0.007806859 &    -7.2199e-3\\ \hline
50  &     0.007817052  &   \ 5.8263e-3\\ \hline
79  &     0.007807948  &   -5.8262e-3\\ \hline
11  &     0.007816705  &   \ 5.3825e-3\\ \hline
118 &     0.007808295 &    -5.3825e-3\\ \hline
115 &     0.007816347 &    \ 4.9250e-3\\ \hline
14  &     0.007808652 &    -4.9249e-3 \\ \hline
\end{tabular}
\end{center}
\caption{Values of 7 MSB's with largest deviations from uniformity
for one orbit of {\tt jsr}{\tt \_}{\tt z}.
}
\label{t:MSBZ}
\end{table}

We note that the {\tt MWC} generator is also not suitable for use with the Ziggurat method.
The reason is that within each orbit, the 7 lsb's of \verb"jsr_w" are not uniformly distributed. Therefore
this PRNG would not choose each of the $128$ intervals at the required uniform probability of 1/128. Indeed numerical experiments
show that the $\chi^2$ statistics on outputs of the Ziggurat method based on this PRNG start
to significantly deviate from the expected value after about $2^{32}$ outputs.

Finally we note that {\tt MWC} is one of the standard uniform random number generators in the
statistical software \verb"R". Given the above non-uniformity of this generator, we caution against its use in simulations.

\section{A statistical analysis of Matlab's randn}

\subsection{The underlying uniform PRNG of randn}
The Matlab software has a built-in function \texttt{randn}
to produce normally distributed random numbers, which is also based on
the Ziggurat method. A matlab code compatible with the pre-compiled built-in function
appears in \cite{matlab_randn}.  In contrast to \cite{Ziggurat}, Matlab's \texttt{randn}
is based on a combination of two different 32-bit registers,
{\tt jsr} and {\tt icng}. Here is a pseudo-C code corresponding
to matlab's \verb"randn":
\begin{verbatim}
unsigned long jsr,icng;
long hz,iz;

#define SHR0 (jsr^=jsr<<13, jsr^= jsr>>17, jsr^=jsr<<5, jsr)
#define CNG  (icng = 69069*icng+1234567)
#define RNOR (hz=CNG+SHR0,iz=hz&63,(fabs(hz)<kn[iz])?hz*wn[iz]:matlab_nfix())
\end{verbatim}
The first register \texttt{jsr} is updated by {\tt SHR0},
as a linear shift register with maximal period
of $2^{32}-1$. The second register \texttt{icng} is updated as a multiplicative
congruential generator, with maximal period $2^{32}$. The output which serves as a uniform random
number is their sum \texttt{(jsr + icng) mod 2$^{32}$}.
We denote the transition function of \verb"jsr" by \verb"T", and that of \verb"icng" by \verb"R". We also denote its
multiplicative part by \verb"R"$_0$, e.g. \verb"R"$_0$\verb"(x)=69069*x mod "$2^{32}$.
Since the periods of the two registers are
relatively prime, the combined period of \texttt{randn} is their product, a number close to
$2^{64}$. Matlab uses a table of size 64, and since it is based on the original Ziggurat publication
\cite{Marsaglia84}, both the points $x_i$, the tables {\tt kn,wn} and the function {\tt matlab\_nfix()} are different
from the ones described in section 3.

Since \verb"icng" is uniformly distributed over $\Omega_{32}$, individual outputs \verb"jsr+icng"
also also uniform over $\Omega_{32}$ and satisfy requirement (\ref{pr_u1}).
However, pairs of consecutive outputs are highly correlated, and fail requirement (\ref{pr_u1_u2}).
Let
$y_1$ and $y_2$ denote two consecutive outputs of this uniform random number generator.
Let $a,b$ denote the unknown initial states at time 1 of the two registers \texttt{jsr} and
\texttt{icng}, that is $(a+b) = y_1$. After a single update of the two registers,
the next output is given by $y_2 = T(a) + R(b)$. However, since $b=y_1-a$, we have that
$y_2 = T(a) - R_0(a) + R(y_1)$. Similar to the analysis of section 3, the transformation $T(a)-R_0(a)$ is highly contractive and
not one-to-one. Therefore, the pair of outputs $(y_1,y_2)$
is not uniformly distributed over $\Omega_{32}\times\Omega_{32}$ and thus fails to satisfy the requirement (\ref{pr_u1_u2}).
Table \ref{t:TaRa} shows the distribution of
$T(a) - R_0(a)$. As shown in the table, some $2^{30.5}$ values are not possible, while other
values are 10 times more probable than expected in a uniform distribution.
Note that $(y_2 - R(y_1))\mbox{ mod } 2^{32}$ is therefore highly
non-smooth, and would not pass a $\chi^2$ test for uniformity.

We now consider the implications of these findings on the resulting normal numbers as computed by {\tt randn}.
Consider, for example, the rejection probabilities at step 3. Since $y$ depends on $x$, the rejection probabilities
deviate slighly from the correct ones. However, when computing these rejection probabilities over large enough
intervals, these $x$-dependent deviations almost cancel out (they are positive for some $x$ and negative for others).
Similarly, tail probabilities at individual $x$-values also deviate from their correct values, but when averaged over large enough intervals
these deviations cancel out. Therefore, even though the underlying generator is not uniform in pairs, its effects
on the resulting normal random numbers is difficult to detect by standard tests.

\begin{table}
\begin{center}
\begin{tabular}{|c|c|}
\hline
\# sources & \# of outputs
\\ \hline
0 & 1590591029 \\ \hline
1 &  1569484236 \\ \hline
2 &  784774346 \\ \hline
3 &  265026908 \\ \hline
4 &  68022535 \\ \hline
5 &  14147755\\ \hline
6 & 2484729\\ \hline
7 & 377496\\ \hline
8 &  51341\\ \hline
9 &  6136\\ \hline
10 & 713\\ \hline
11 & 65\\ \hline
12 & 6\\ \hline
13 & 1\\ \hline
\end{tabular}
\end{center}
\caption{Distribution of the number of sources of $Tx-R_0x$.}
    \label{t:TaRa}
\end{table}

\subsection{Initialization Issues}

We now consider the initialization of {\tt randn} and its possible consequences. Matlab provides two different initialization options,
\[
{\tt randn('state',a);\ \ \ \ \ \ \ \        OR\ \ \ \ \ \ \ \  randn('state',[a\ b]'); }
\]
The first sets the initial value of {\tt jsr} to \verb"a", with the initial value of {\tt icng} set to
a fixed value 362436069. The second option allows to set also the initial value of {\tt icng} to \verb"b".

In many applications, such as parallel computations and stochastic simulations,
there is a need to create many independent sequences of normal random numbers.
In the case of parallel computer systems, it is quite common to initialize the seed of processor
number {\tt id} with a seed of the form {\tt seed$_0$ + id}.
Quite a few works describe the dangers and possible pitfalls in using
sequences of random numbers produced by different initializations of the same generator (see \cite{Brent,Hellekalek,{Pagnutti88}} and
references therein). We now present
a simple example of such a pitfall for matlab's {\tt randn}.

Suppose we wish to simulate $2^{16}$ different paths of a stochastic system that requires normal
random numbers on a parallel computer with 256 processors. A possible code can be for example
\begin{verbatim}
for i=1:256
    for j=1:256
        send to processor i the following:
                    randn('state',[i j]');
                    simulate_random_path();
    end
end
\end{verbatim}
This code ensures that each simulation thread obtains a different seed. However, consider
the output sequences resulting from two initializations {\tt [i j]} and {\tt[i j+64]}. These two initializations
have the same initial value for {\tt jsr} and differ only by the initialization of the register {\tt icng}. Due to the structure of the transition
function \verb"R" of this register, it follows that for all subsequent times, both of these sequences
will have the same six low significant bits. Therefore, neglecting the possible misses in the
Ziggurat method, which require a call to \verb"matlab_nfix()", both sequences will choose the
same indexes (!), and the resulting normal numbers will be highly correlated.

These correlations between output sequences initialized with different but related seeds are due to the failure of {\tt randn} to
satisfy requirement (\ref{pr_u_utag}). However, we might be tempted to conclude that if instead we perform initializations
with random unrelated seeds the resulting sequences will be uncorrelated. However, even in such
cases there can be non-trivial correlations between the first few values of different output sequences.
Consider the output sequences from $n$ initializations of the form
\verb"randn('state',v[i])", where $\{v_i\}_{i=1}^n$ is a set of $n$ random non-zero 32 bit integers.
Suppose there are four distinct indices, $i,j,k,l$, such that $v_i\oplus v_j = v_k \oplus v_l = \alpha$,
where $\oplus$ denotes bitwise exclusive or. Then,
for the first few outputs, the resulting 4-tuples of outputs are not independent.
To see this, denote $x^i,x^j,x^k,x^l$ the resulting first output of the uniform random number generator
initialized with $v_i,v_j,v_k,v_l$, respectively. Since all these initializations have the same initial
value for the register {\tt icng}, the first output is given by
\[
x^i = T(v_i) + R(icng)\quad x^j = T(v_j) + R(icng) = T(v_i)\oplus T(\alpha) + R(icng)\\
\]
with similar expressions for $x^k$ and $x^l$. As an example, assume that the most significant bit of
$T(\alpha)$ is zero. Then, with high probability the most significant bit of $x^i$ and $x^j$ (and of $x^k$ and $x^l$)
will be the same. Therefore, the resulting normal numbers will have the same sign. For the second output the sign
will be determined by the most significant bit of $T^2(\alpha)$, etc. (assuming that no intermediate calls to
\verb"matlab_nfix()" occurred). Another example of dependency occurs if we consider the 7 least significant bits of
$T(\alpha)$. For simplicity, if these are all zeros, then the numbers $x^i$ and $x^j$ (and $x^k$ and $x^l$)
point to the same index for the Ziggurat method. This again leads to correlations between the 4-tuples
of normal outputs. Needless to say, such correlations between different streams may bias a stochastic simulation
in unexpected ways. Given $n$ initializations, the average number of such 4-tuples is of the order of
${n \choose 4}2^{-32}$. Therefore, after only $n=O(512)$ random sequences there will be on average one such 4-tuple.

\section{Summary and Discussion}

In this paper we presented various statistical weaknesses in the
published implementation of the Ziggurat and Monty Python methods and in the underlying uniform PRNG
of matlab's built-in function \texttt{randn}.
As also noted in other works, the main take home lessons from our analysis are:
i) The set ${\cal S}$ of internal states of the generator should be much larger than the output set $\Omega$.
ii) Correlations between consecutive outputs of a uniform RNG can have detrimental
effects on the results of a stochastic simulation.
iii) The creation of many different sequences of random numbers via initializations with different seeds
must be done with great care.

Regarding the initialization of random number generators, we note that most implementations
use $P = I$ or some other relatively simple scheme, in which the seed is typically entered into the inner state
in a linear fashion. However, since initialization is done only rarely
it is possible to spend many more CPU cycles on this stage, and make the inner state be dependent on the initial seed
in a much more complicated and non-linear manner.

We remark that there is an interesting connection between our analysis of PRNG's and cryptanalysis of
stream and block cyphers. For example, our analysis of the statistical correlations of PRNG's
initialized with different but related seeds is similar to the 'related key attacks' introduced by Biham in \cite{Biham93}, and used to crack
the WEP wireless encryption protocol \cite{WEP}. This serves as yet another
justification for making the initialization stage (the function $P$ in the notation of section 2)
a complicated non-linear function.

There is yet another connection between PRNG's and cyphers. In cryptography,
the designer would like the security of a cypher system not be
dependent on its specific initialization by the user (e.g., with say counters or
initial values (IV's)
increasing by one). We submit that a similar requirement should hold
in the design of a PRNG. The output of a PRNG (and correlations between different
output runs) should not be highly dependent on simple and natural initializations by the user,
who typically does not know nor wishes to fully understand the inner workings of
the random number generator at his disposal.

\noindent
{\bf Acknowledgments:} The author would like to thank Prof. Adi Shamir for interesting discussions.

\section*{Appendix}

Let $X$ and $Y$ be two discrete random variable over $k$ possible values, with probability distributions
$(p_1,\ldots,p_k)$ and
$(q_1,\ldots,q_k)$, respectively. Our aim is to estimate the number of outputs needed from i.i.d. realizations
of $Y$ to check the hypothesis if $Y$ has the same probability distribution as $X$, using
the $\chi^2$ test with $k$ bins. Let $(y_1,\ldots,y_N)$ be $N$ random samples from
the distribution of $Y$, and let $(z_1,\ldots,z_k)$ denote the number of occurrences of the values
$(1,\ldots,k)$ in the sequence $\{y_i\}_{i=1}^N$. Then the $\chi^2$ statistic is given by
\[
T = \sum_{i=1}^k \frac{(z_i-Np_i)^2}{Np_i} = \sum_{i=1}^k \frac{z_i^2}{Np_i} - N
\]
Its mean (expected) value is
\begin{equation}
  \mathbb{E}T = \sum_{i=1}^k\frac{\mathbb{E}z_i^2}{N p_i} - N
\end{equation}
If $Y\sim(q_1,\ldots,q_k)$ then each $z_i$ follows a Binomial distribution $Bin(N,q_i)$. Therefore,
\begin{equation}
  \mathbb{E}T = \sum_{i=1}^k \frac{N^2 q_i^2 + N q_i (1-q_i)}{N p_i} -N
\end{equation}
Writing $q_i = p_i(1 + \ve_i)$ gives
\begin{equation}
  \mathbb{E}T = (k-1) + (N-1)\sum_{i=1}^k p_i \ve_i^2 + \sum_{i=1}^k \ve_i
\end{equation}
Since a $\chi^2$ distribution with $k$ degrees of freedom has a variance of $2k$,
it follows that to distinguish between the distribution of $X$ and $Y$, the $\chi^2$ statistic
must significantly deviate from $k-1 + \sqrt{2k}$. Thus, we require that
\begin{equation}
  N = O\left( \frac{\sqrt{2k}}{\sum_i p_i \ve_i^2} \right)
\end{equation}

\end{document}